\documentclass{elsart}

\usepackage{amssymb}
\usepackage{amsfonts}
\usepackage{mathrsfs}  
\usepackage{amsmath}

\newtheorem{theorem}{Theorem} 

\newtheorem{proposition}[theorem]{Proposition}
\newtheorem{corollary}[theorem]{Corollary}

\newtheorem{definitions}[theorem]{Definitions}

\newenvironment{proof}{\noindent{\em Proof~}}{\hfill $\Box$}

\newtheorem{preexample}[theorem]{Example} 
\newenvironment{example}{\begin{preexample}\em}{\end{preexample}}

\makeatletter
\newcommand{\mysubsection}%
{\@startsection{subsection}{2}{\z@}{-3.25ex plus -1ex minus -.2ex}{-1ex}{\normalsize\sc}}
\makeatother

\renewcommand{\phi}{\varphi}
\renewcommand{\epsilon}{\varepsilon}
\newcommand{\U}{\mathbb{U}}

\newcommand{\bdu}{\partial \U}
\newcommand{\bilin}[2]{\langle#1,\,#2\rangle}

\newcommand{\C}{\mathbb{C}}

\newcommand{\chat}{\hat{\C}}

\newcommand{\cphi}{C_\phi}

\newcommand{\cphistar}{\cphi^*}

\newcommand{\eqdef}{:=}

\newcommand{\fhatn}{\hat{f}(n)}
\newcommand{\fhat}{\hat{f}}

\newcommand{\ghatn}{\hat{g}(n)}
\newcommand{\goesto}{\rightarrow}

\newcommand{\htwo}{H^2}

\newcommand{\norm}[1]{\|#1\|}

\newcommand{\omegaphi}{\Omega_\phi}

\newcommand{\phibar}{\overline{\phi}}
\newcommand{\phinfty}{\phi(\infty)}
\newcommand{\phinftybar}{\overline{\phinfty}}

\newcommand{\phie}{\phi_e}

\newcommand{\phinv}{\phi^{-1}}
\newcommand{\phieinv}{{\phie}^{-1}}

\newcommand{\phieprime}{{\phie}'}

\newcommand{\ratu}{{\rm Rat}(\U)}

\newcommand{\reg}{{\rm reg}}

\newcommand{\Rinv}{R^{-1}}

\newcommand{\rhoz}{\rho(z)}

\newcommand{\rhowj}{\rho(w_j)}

\newcommand{\seq}[1]{\{#1\}}
\newcommand{\sigmaseq}{\{\sigma_j\}_{j=1}^d}

\newcommand{\ubar}{\overline{\U}}
\newcommand{\Ue}{\U_e}

\newcommand{\wbar}{\bar{w}}
\newcommand{\wjbar}{\overline{w_j}}
\newcommand{\wjbarinv}{1/\wjbar}

\newcommand{\zbar}{\overline{z}}

\begin{document}

\begin{frontmatter}

\title{Adjoints of rationally induced \\composition operators}

\author{Paul S. Bourdon}
\address{Department of Mathematics\\  Washington and Lee University, Lexington VA 24450}
\ead{pbourdon@wlu.edu}

\author{Joel H. Shapiro}

\address{Department of Mathematics and Statistics\\ Portland State University, Portland, OR 97207}
\ead{shapiroj@pdx.edu}

\begin{abstract}
We give an elementary proof of a  formula recently obtained by Hammond, Moorhouse, and Robbins  for the adjoint of a rationally induced composition operator on the Hardy space $\htwo$ [Christopher Hammond, Jennifer Moorhouse, and Marian E. Robbins, {\em Adjoints of composition operators with rational symbol,} J. Math. Anal. App., 341 (2008) 626--639]. We discuss some variants and implications of this formula, and use it to provide a sufficient condition for a rationally induced composition operator adjoint to be  a compact perturbation of a weighted composition operator.
\end{abstract}

\end{frontmatter}


\section*{Introduction}\label{introduction}
We study composition operators $\cphi:f\goesto f\circ\phi$ acting on the Hardy space $\htwo$ of the open unit disc $\U$ of the complex plane. Our goal is to provide a simple derivation for an intriguing formula, established recently by Hammond, Moorhouse, and Robbins \cite{HMR}, that describes the adjoints of composition operators induced by rational selfmaps of $\U$.  

The study of composition operator adjoints was initiated more than twenty years ago by Carl Cowen \cite{C}, who showed that if $\phi$ is linear-fractional then $\cphistar$, the adjoint of $\cphi$ on $\htwo$,  has the form $M_gC_\sigma M_h^*$, where $M_g$ and $M_h$ are the operators of multiplication by simple rational functions $g$ and $h$, bounded on $\U$, and $\sigma$ is a linear-fractional selfmap of $\U$---all these functions being easily obtainable from the coefficients of $\phi$ (see also \cite[Theorem 9.2, page 322]{CMB}).

Cowen's work showed, in particular,  that one should expect the adjoint of a composition operator to be substantially more complicated than just another composition operator; in the linear-fractional case it is the adjoint of a multiplication operator followed by a {\em weighted composition operator\/}  $M_gC_\sigma$. Complications notwithstanding, Cowen's formula has provided an essential tool in a number of investigations involving linear-fractionally induced composition operators, for example:  subnormality and co-subnormality \cite{C}, computation of norms (see Hammond's fundamental paper \cite{H}, and subsequent work \cite{ BR, BHFS, EJRS} based upon it), essential normality \cite{BLNS}, self-commutator properties \cite{BM}, $C^*$ algebras \cite{KMM1,KMM2}, and singular value decompositions \cite{CD}. 

Useful as Cowen's work has been,  the study of adjoint formulas for composition operators has only recently been pushed beyond the linear-fractional case.  Although several different approaches have been tried (see \cite{MV} and the references described therein), explicit formulas have been hard to come by, even for rationally induced operators.  

 Recent work of Cowen and Gallardo \cite{CG} (see also the discussion of  \cite{CG} in \cite{HMR}) made the point that for general rational selfmaps of the disc the determination of $\cphistar$ must involve what they called ``multiple-valued weighted composition operators.'' For a graphic demonstration of how multiple-valuedness can arise, consider the simplest rational selfmap of $\U$ that is not linear-fractional: $\phi(z) = z^2$. Elementary calculations show that for each function $f(z) = \sum_{n=0}^\infty\fhat(n)z^n$ in $\htwo$, 
$$
    \cphistar f(z) ~=~ \sum_{n=0}^\infty \fhat(2n)z^n 
                       ~=~  \frac{1}{2}[f(\sqrt{z}) + f(-\sqrt{z})]  \qquad (z\in\U).
$$
Thus for this example we can regard $\cphistar$  as a formal average of  ``composition operators'' induced by branches of the square root function, with the inducing maps having singularities along a branch cut which, even when subjected to composition with an arbitrary $f\in\htwo$, retain enough cancellation to yield a result that is holomorphic on the whole unit disc. 

Hammond, Moorhouse, and Robbins \cite{HMR} recently corrected a formula in \cite{CG} for the adjoint of a composition operator on $H^2$ induced by a rational function $\phi$. The method of both papers  centers on Cauchy integral representations of the $H^2$ inner product $\langle f, \cphistar g\rangle$, for polynomials $f$ and $g$.  A case-based application of Cauchy's Theorem followed by a change of variables results in a representation of the inner product involving a generalization of the map $\sigma$ that figured in Cowen's  linear-fractional theorem, the difficulty being that beyond the linear-fractional case this generalization turns out to be {\em multiple-valued,} so that its branching behavior needs to be taken into account.

Our initial goal here is to obtain the Hammond-Moorhouse-Robbins (HMR) formula in a straightforward algebraic fashion, resulting in a unified derivation that is considerable shorter than that in \cite{HMR}. After setting out the necessary background material in the next section, we devote \S \ref{maintheorem_section} to our proof of the HMR  formula,  discussing some of its variants in \S \ref{variants}. In \S\ref{regularity} we address the question of when $\cphistar$ is a  sum involving ``legitimate'' weighted composition operators,  after which we apply our results to the question of which rationally induced composition operators have adjoint equivalent, modulo compacts, to a weighted composition operator. 
	
\section{Background material} \label{background}
For expositional completeness we collect here some standard background material that will be needed in what follows.

 \mysubsection{Mapping rationally.}\label{ratselfmap_intro}
 We will need to think of rational functions, even self-maps of the unit disc, as mappings taking the extended complex plane (i.e., the Riemann sphere) $\chat$ into itself. Throughout this paper, $\phi$ denotes an element of $\ratu$, the collection of all rational functions of one complex variable that map $\U$ into itself.  Associated with $\phi$ is its ``exterior'' map $\phie \eqdef \rho\circ\phi\circ\rho$, where $\rho:\chat\goesto\chat$ denotes inversion in the unit circle ($\rho(z) \eqdef 1/\zbar$). Thus $\phie$ maps the exterior $\Ue$ of the closed unit disc (viewed as a subset of $\chat$) into itself. On the other hand it is easy to check that $\phieinv(\U)\subset \U$, a fact that will be of utmost importance to us. 

\mysubsection{Regular and critical values.}\label{reg_crit}
For any rational function $R$, expressed as a quotient of relatively prime polynomials, the {\em degree\/} of $R$ is the larger of the degrees of its numerator and denominator. It is an elementary exercise  to show that if the degree of $R$ is $d$ then for each point $w\in\chat$ the inverse image $R^{-1}(\{w\})$ has, counting multiplicities, exactly $d$ points (see \cite[Problems 25--32, pp. 181--182]{F} for example).  If $R^{-1}(\{w\})$ has $d$ {\em distinct\/} points we will say $w$ is a {\em regular value\/} of $R$. Clearly, for any rational function, all but finitely many points of $\chat$ are regular values. We call those  values that are not regular {\em critical\/} values (perhaps more commonly these are called ``branch points of $R^{-1}$).

By elementary function theory,  $w\in\chat$ is a critical value of $R$ iff $R^{-1}(\{w\})$ contains at least one point having no neighborhood on which $R$ is univalent. Such a point is called a {\em critical point\/}. Thus a value is critical for $R$ if and only if its preimage contains a critical point. 

Suppose that $w$ is a regular value of $R$, so that $R^{-1}(\{w\})$  consists of distinct points $\seq{z_1, z_2, \, \ldots \, z_d}$ in $\chat$. Then, again by elementary function theory,  each point $z_j$ is the center of an open disc $\Delta_j$ to which $R$ has a univalent restriction, hence the restriction of $R$ to $\Delta_j$ has a holomorphic inverse $\sigma_j$ on $R(\Delta_j)$; it is called a  {\em branch\/} of  $\Rinv$ defined on $R(\Delta_j)$. Upon taking the intersection of the $R$-images of the $d$ discs $\Delta_j$ we obtain a common neighborhood of $w$ on which all the distinct branches $\sigmaseq$ are defined.

\mysubsection{Branching out.}\label{branching_out}
  By way of introduction, suppose more generally that $V$ and $W$ are arbitrary sets and $f:V\goesto W$ is a function that is one-to-one on a subset $V_1$ of $V$. Then the restriction of $f$ to $V_1$ has an inverse $\sigma_1$ mapping $\Omega:=f(V_1)$ back onto $V_1$. Let's call $\sigma_1$ a {\em branch of $f^{-1}$ defined on $\Omega$.} It's easy to check that if $\sigma_2$ is another branch of $f^{-1}$ defined on $\Omega$, with $V_2=\sigma_2(\Omega)$, then $\sigma_1\equiv\sigma_2$ on $f(V_1\cap V_2)$. 

Suppose, in particular that $V$ is an open subset of $\chat$,  $f$ is holomorphic on $V$, and that $\sigma_1$ and $\sigma_2$ are branches of $f^{-1}$ defined on a connected open subset $\Omega$ of $\chat$. Then $\sigma_1$ and $\sigma_2$ are both holomorphic on $\Omega$, the sets $V_1$ and $V_2$ of the last paragraph being, respectively, $\sigma_1(\Omega)$ and $\sigma_2(\Omega)$. Thus if $V_1\cap V_2$ is nonempty, then its image under $f$ is, by the last paragraph,  a nonvoid open subset of $\Omega$ on which $\sigma_1\equiv\sigma_2$. By the uniqueness theorem for holomorphic functions, $\sigma_1\equiv\sigma_2$ on all of $\Omega$. In summary:

\begin{proposition} \label{disjointness_prop}
Suppose $f$ is holomorphic on an open subset of $\chat$ and $\sigma_1, \sigma_2$ are two branches of $f^{-1}$ defined on a connected open subset $\Omega$ of $\chat$. Then either  $\sigma_1(\Omega)\cap\sigma_2(\Omega)=\emptyset$  or $\sigma_1\equiv \sigma_2$ on $\Omega$.
\end{proposition}

\mysubsection{Continuing branches.}\label{continuing_branches}
Suppose $R$ is a rational function. Let $\reg(R)$ denote the collection of points in $\chat$ that are regular values of $R$. The discussion of regular values in \S\ref{reg_crit} shows that $\reg(R)$ is an open subset of $\chat$; furthermore $\reg(R)$ supports ``arbitrary continuation''  of branches of $\Rinv$. More precisely:
\begin{proposition}  \label{continuation_prop}
Suppose $R$ is a rational function and $\gamma:[0,1]\goesto\reg(R)$ a continuous curve. Suppose $\sigma$ is a branch of $\Rinv$ defined on an open disc centered at  $\gamma(0)$ and lying in $\reg(R)$. Then $\sigma$ has a holomorphic continuation along $\gamma$.
\end{proposition}
   \begin{proof}
     Forget the curve for a moment and just suppose we have two open discs $\Delta_1$ and $\Delta_2$ contained in $\reg(R)$, each of which supports $d$ distinct branches of $\Rinv$; say $\seq{\sigma_j^{(1)}}_{j=1}^d$ on $\Delta_1$ and $\seq{\sigma_j^{(2)}}_{j=1}^d$ on $\Delta_2$. Suppose $\Delta_1\cap\Delta_2\neq\emptyset$.
  
{\em Claim:} For each index $j$ the branch $\sigma_j^{(1)}$ of $\Rinv$ has a holomorphic extension to $\Delta_1\cup\Delta_2$.

{\em Proof of Claim.} Fix $w_0\in \Delta_1\cap\Delta_2$. Then for $i=1,2$;
     $$
         \Rinv(\seq{w_0}) = \seq{\sigma_j^{(i)}(w_0)}_{j=1}^d~,
     $$
 where the set on the right has $d$ distinct elements. Thus there is an index $k$ such that $\sigma_j^{(1)}(w_0)=\sigma_k^{(2)}(w_0)$.    Now both $\sigma_j^{(1)}$ and $\sigma_k^{(2)}$ are branches of $\Rinv$ defined on the connected  open set $\Delta_1\cap\Delta_2$, and by our choice of $k$, the images of that intersection under the two branches has nonvoid intersection. Thus by Proposition \ref{disjointness_prop}, $\sigma_j^{(1)}\equiv \sigma_k^{(2)}$ on the intersection, and so the function $\sigma$ defined on $\Delta_1\cup\Delta_2$ by: 
 $$
     \sigma ~= ~ \left\{
     \begin{array}{ccc}
            \sigma_j^{(1)}&\mbox{on}&\Delta_1\\
            \sigma_k^{(2)}&\mbox{on}& \Delta_2   
      \end{array}
      \right.                 
 $$
 gives the desired holomorphic extension of $\sigma_j^{(1)}$ to $\Delta_1\cup\Delta_2$.
 
 Now let's talk about the the curve $\gamma$. By its compactness in the open set $\reg(R)$ and the discussion of \S\ref{reg_crit}, it can be covered by a chain of open discs $\seq{\Delta_j}_{j=1}^n$ lying in $\reg(R)$, each supporting $d$ distinct branches of $\Rinv$, with  $\Delta_0$ being the open disc on which the branch $\sigma$ of $\Rinv$ is initially defined, and each disc having nontrivial intersection with its successor. Then successive application of the {\em Claim} effects a holomorphic continuation of $\sigma$ along the curve $\gamma$.
 \end{proof}

 \begin{corollary}  \label{sc_continuation}
  Suppose $R$ is a rational function and $\sigma$ is a branch of $\Rinv$ defined on some disc $\Delta\subset\reg(R)$. Suppose $W$ is a {\em simply connected\/} domain containing $\Delta$ and lying in $\reg(R)$. Then $\sigma$ has a holomorphic continuation to a branch of $\Rinv$ defined on $W$. 
\end{corollary}
\begin{proof}
Proposition \ref{continuation_prop}, along with the Monodromy Theorem \cite[Theorem 10.3.1, page 12]{Hi}, yields the desired continuation of $\sigma$ to $W$. That this continuation is a branch of $\Rinv$ follows quickly from the uniqueness theorem for holomorphic functions (or, perhaps more ostentatiously, from the ``law of permanence of functional equations'' \cite[\S10.7, pp. 31--36]{Hi}). 
\end{proof}

\mysubsection{$\htwo$ and its reproducing kernels.}\label{repro_kernels}
The space $\htwo$ is the collection of functions holomorphic on $\U$ with square-summable Maclaurin-coefficient sequence. The norm $\norm{f}$ of a function $f\in\htwo$ with Maclaurin series $f(z) = \sum_{n=0}^\infty \fhatn z^n$ is defined by:
$$
    \norm{f}^2= \sum_{n=0}^\infty |\fhatn|^2~;
$$
this norm makes $\htwo$ into a Hilbert space with inner product
$$
     \bilin{f}{g} = \sum_{n=0}^\infty \fhat(n)\overline{\ghatn}
         \qquad (f,g\in\htwo).
$$ 
A simple computation shows that for each $z\in\U$ and $f\in\htwo$,
\begin{equation} \label{repro_ker_eqn}
      f(z) = \bilin{f}{K_z} 
               \qquad \mbox{where}
                      \qquad K_z(w) \eqdef \frac{1}{1-\zbar w}  \qquad (w\in\U).                   
\end{equation}
The function $K_z$ (clearly in $\htwo$)  is called the {\em reproducing kernel} for  the point $z$.

\mysubsection{The ``kernel'' of our argument.}\label{kernel_of_arg}
The importance to our work of reproducing kernels  stems from the observation that
\begin{equation}\label{crucial_observation}
    \cphistar f(z) = \bilin{\cphistar f}{K_z} = \bilin{f}{\cphi K_z}.
\end{equation}
Now $\cphi K_z(w) = (1-\zbar\phi(w))^{-1}$ is a rational function of $w$ whose partial fraction expansion we will determine explicitly. The desired formula for $\cphistar f(z)$ will emerge upon substitution of this expansion into \eqref{crucial_observation}; the details are in the next section.

\mysubsection{Much ado about zero.}\label{much_ado}
Although the general formula for $\cphistar f(z)$ can be quite complicated, for $z=0$ it could not be simpler. Upon substituting $z=0$ into \eqref{crucial_observation} and noting that $K_0\equiv 1$, we obtain
$$
    \cphistar f(0) = f(0) \qquad (f\in\htwo),
$$
an observation we will need several times in the sequel.

\mysubsection{The backward shift.} \label{backward_shift}
In \S\ref{variants} we will derive alternate formulae for $\cphistar f$, one of which will involve the backward shift operator $B$, defined for $f\in\htwo$ by:
$$
          (Bf)(z) \eqdef 
            \left\{
          \begin{array}{cll}
          {\displaystyle\frac{f(z)-f(0)}{z}}
                      & \quad \mbox{if}\quad& 
                                  z\in\U\backslash\{0\} \\
                     ~&              ~              &~                                   \\     
          f'(0)      & \quad \mbox{if}\quad& z=0.
          \end{array}
          \right.
$$
This operator is so named because it has the effect of shifting Maclaurin coefficient sequences ``back one step,'' i.e.,
$$
     B:\sum_{k=0}^\infty \fhat(k)z^k \goesto \sum_{k=0}^\infty \fhat(k+1)z^k
     \qquad (f\in\htwo).
$$
It is easy to check that $B$ is the $\htwo$-adjoint of the multiplication operator $M_z$ defined, albeit with some abuse of notation, by:
  $$
      M_z: f(z)\goesto zf(z) \qquad (f\in\htwo).
  $$
Because $M_z$ has the effect of moving Maclaurin coefficient sequences ``forward one step'' it is called---not surprisingly---the ``forward shift.''

\section{Main Theorem}  \label{maintheorem_section} Our goal in this section is to give an elementary proof of the formula for $\cphistar$ proved by Hammond, Moorhouse, and Robbins \cite{HMR}. 

\mysubsection{The exterior map.} \label{exterior_map} In order to state the result efficiently, recall from \S\ref{ratselfmap_intro} that each rational selfmap $\phi$ of $\U$ has a companion rational selfmap $\phie=\rho\circ\phi\circ\rho$ of $\Ue$, where $\rho$ is the mapping of inversion in the unit circle. Note that the degree of $\phie$ is the same as that of $\phi$, and that a point of $\chat$ is a regular value of $\phie$ if and only if its reflection in the unit circle is a regular value of $\phi$. 

\mysubsection{The Main Theorem.}\label{main_thm_statement} 
Recall that if $z_0\in\U$ is a regular value of $\phie$ and $\phie$ has degree $d$, then, as we saw in \S\ref{reg_crit}, there is a neighborhood  $V$ of $z_0$, which we assume to be connected, on which are defined $d$ distinct branches $\sigmaseq$ of $\phieinv$.  The univalence of these branches insures that their derivatives vanish nowhere on $V$.  Note that Proposition~\ref{disjointness_prop} shows each element of $V$ is a regular value of $\phie$ and that our earlier observation  $\phie^{-1}(\U)\subset \U$ shows  $\sigma_j(z)\in \U$ for every $z\in V\cap\U$.

\begin{theorem}[\cite{HMR}] \label{maintheorem}
Let $\phi\in\ratu$ have degree $d>0$ and let $f\in\htwo$. Suppose $z_0\in\U$ is a regular value of $\phie$ and  $V\subset \U$ is any connected neighborhood of $z_0$ on which are defined $d$ distinct branches $\sigmaseq$ of $\phieinv$.  

{\em (a)} If $\phinfty\neq\infty$, then for all $z\in V\setminus\{ 1/\phinftybar\}$
\begin{subequations}\label{HMR_formulas}
  \begin{equation}\label{HMR_formula}
        \cphistar f(z) = \frac{f(0)}{1-\phinftybar \, z} ~+~ 
                z \sum_{j=1}^d \frac{\sigma_j'(z)}{\sigma_j(z)} \, f(\sigma_j(z))~.
   \end{equation}
{\em (b)} If $\phinfty=\infty$, then for all $z\in V$ 
    \begin{equation} \label{HMR_formula_infty}
        \cphistar f(z) ~=~ \left\{ 
             \begin{array}{cll}
               {\displaystyle z \sum_{j=1}^d 
                     \frac{\sigma_j'(z)}{\sigma_j(z)} \, f(\sigma_j(z))}
                     &\quad \mbox{if} \quad & z\neq 0 \\
                     & & \\
                f(0) & \quad \mbox{if}  \quad  & z=0  
                               \end{array} \right.
    \end{equation}
    \end{subequations}
\end{theorem}

\mysubsection{Discussion.}\label{discussion}
   Before we begin a rigorous proof, some words need to be said about these formulas. Note first that they must reflect, as they do, the fact that $\cphistar f(0) = f(0)$ for every $f\in\htwo$ (see \S\ref{much_ado}).  In  formula (\ref{HMR_formula_infty}), since $\cphistar f$ is continuous at $0$, the second line on the right must be a limiting case of the first one.  This may be proved independently as follows.
      
{\em Formula  {\em (\ref{HMR_formula_infty})} and continuity at the origin.} The first  line of (\ref{HMR_formula_infty}) appears as a natural way of interpreting (\ref{HMR_formula}) when $\phi(\infty)=\infty$. In this case $\phie(0)=0$ so, assuming that $V$ contains $0$, we conclude that  exactly one of the $\sigma_j$'s has value zero at the origin (Proposition~\ref{disjointness_prop}). For this $j$ we see that  $z\sigma_j'(z)/\sigma_j(z)\goesto \sigma_j'(0)/\sigma_j'(0) =1$ as $z\goesto 0$ (recall that univalence of $\sigma_j$ ensures $\sigma_j'(0)\ne 0$). The other such terms  clearly converge to $0$ as $z\goesto 0$, hence the first formula of (\ref{HMR_formula_infty}) does indeed tend to $f(0)$ as $z\goesto 0$.

{\em Remarks on equation {\em (\ref{HMR_formula})}.} Since $\phie(0) = 1/\phinftybar$,  the hypothesis $\phinfty\neq\infty$ means that $\phie(0)\neq 0$ while the hypothesis $z\neq 1/\phinftybar$ means that no denominator $\sigma_j(z)$ on the right-hand side of (\ref{HMR_formula}) is zero. If, however, $\phie(0)$ were to lie in $V$,  and we were to allow $z$ to take that value, then one (and exactly one, by Proposition~\ref{disjointness_prop}) of the $\sigma_j(z)$'s would be zero. 
    Since $\sigma_j'(z)\neq 0$ this would, since $z \neq 0$,  introduce in the second term on the right-hand side of (\ref{HMR_formula}) a pole of order 1 at $z$. But now $\phi(\infty) = 1/\zbar\in\Ue$, which provides the first term on the right-hand side of (\ref{HMR_formula}) with a pole of order 1 at $z$, and since the left-hand side of that equation is holomorphic on $\U$, these poles must cancel each other.  Thus formula (\ref{HMR_formula}) could be considered valid at $z=1/\phinftybar$ in the sense that the left-hand side is the limit of the right as $z$ approaches $1/\phinftybar$. 
    
{\em A convenient choice of $V$.}  Let $\widetilde{V}$ be the unit disk with radial slits from each critical value of $\phie$ to the unit circle removed.  Then $\widetilde{V}$ is a simply connected domain and Corollary~\ref{sc_continuation} shows that $\phie^{-1}$ has $d$ distinct branches $\{\sigma_j\}_{j=1}^d$ defined on $\widetilde{V}$. We can assert that formula (\ref{HMR_formula}) is valid for all $z\in \widetilde{V}$ provided we interpret the formula according to the discussion provided above.   Note that if $\phie$ has no critical values in $\U$, then  $\widetilde{V} = \U$ and the branches $\{\sigma_j\}_{j=1}^d$ of $\phie^{-1}$ become holomorphic self-maps of $\U$.  

\mysubsection{Proof of Main Theorem.}\label{proof_main_thm}  Let $z_0\in\U$ be a regular value of $\phie$ and  $V\subset \U$ be any connected neighborhood of $z_0$ on which are defined $d$ distinct branches $\sigmaseq$ of $\phieinv$.  Let $z\in V$ so that $z$ is a regular value of $\phie$ whose $d$ preimages under $\phie$ are, of course, $\{\sigma_j(z)\}_{j=1}^d$.  Given the discussion of \S\ref{discussion},  we may assume that $z\ne 0$ and we shall also assume that  that $z\neq 1/\phinftybar$, so that $\cphi K_z(w)$ has a finite limit as $w\goesto\infty$.

Since $z\neq 0$ it is easy to check that $\cphi K_z(w)$, viewed as a rational function of $w$, has the same degree as $\phi$ (and therefore of $\phie$); it has its poles at the points $w$ for which $\phi(w)=1/\zbar\eqdef\rho(z)$, i.e. at the points of $\phinv(\seq{\rhoz})$. Note these poles are simple since $1/\zbar$ is a regular value of $\phi$.  As we pointed out in \S\ref{kernel_of_arg}, the crucial step of our proof will be to determine the partial fraction expansion of $\cphi K_z$.

Because  $z\in\U$ is a non-zero regular value of $\phie$, its reflection $\rhoz\in\Ue$ is a finite regular value of $\phi$, with $\phinv(\rhoz)$ consisting of the $d$ distinct points in $\{1/\overline{\sigma_j(z)}\}_{j=1}^d$, all of which are finite by our  assumption that  $z\neq 1/\phinftybar$.  For notational convenience, we set 
$$
1/\overline{\sigma_j(z)} = w_j
$$
for each $j\in \{1, 2, \ldots, d\}$. 
 
 Thus we have this partial fraction expansion for $\cphi K_z(w) = 1/(1-\zbar\,\phi(w))$:   
\begin{equation}\label{first_pfe}
   \cphi K_z(w) ~=~  \alpha ~+~\sum_{j=1}^d \frac{\beta_j}{w-w_j}
                        ~=~ \alpha ~+~ \sum_{j=1}^d \frac{\beta_j/w_j}{w/w_j-1} ~ ,                       
\end{equation}
where 
\begin{equation} \label{alpha_calcn}
     \alpha ~=~ K_z(\phi(\infty)) ~=~ \frac{1}{1-\zbar\phi(\infty)} \in\C,
\end{equation}
with the $\beta_j$'s, which also depend on $z$, to be determined shortly. Note for further reference that $\alpha=0$ in case $\phi(\infty)=\infty$. 

Since each of the points $\rhowj$ is in $\U$, equation (\ref{first_pfe}) can be rewritten in terms of the reproducing kernels for these points:
\begin{equation}\label{reproker_eqn}
    \cphi K_z(w) 
        ~=~ \alpha  ~-~ \sum_{j=1}^d \frac{\beta_j}{w_j}K_{\rhowj}(w) ~ .
\end{equation}
Pending determination of the coefficients $\beta_j$, the following calculation will produce the formula for $\cphistar$. First, substitute (\ref{reproker_eqn}) into  (\ref{crucial_observation}) to obtain, for $f\in \htwo$,
$$
  \cphistar f(z) 
        = \bilin{f}{\cphi K_z} = \bilin{f}{\alpha  - \sum_{j=1}^d \frac{\beta_j}{w_j}K_{\rhowj}}~.
$$
On the right-hand side of this equation interchange the sum and inner product, and use the reproducing kernel formula (\ref{repro_ker_eqn}) to evaluate the resulting inner products.  The result, in view of  (\ref{alpha_calcn}), is:
\begin{equation} \label{second_cphistar_formula}
   \cphistar f(z) 
          ~=~\frac{f(0)}{1-\overline{\phi(\infty)}z} ~-~ 
                \sum_{j=1}^d   \, \overline{\beta_j} \, \rhowj  \, f(\rhowj) ~ .
\end{equation}

\bigskip\noindent
{\em Evaluation of the $\beta_j$'s.}
Because we are assuming that the poles $\seq{w_j}_{j=1}^d$ of $K_z\circ\phi$ are all simple, this is an easy residue calculation. From  (\ref{first_pfe}) it's just:
\begin{equation}\label{rescalcn}
   \beta_j = \lim_{w\goesto w_j} \frac{(w-w_j)}{1-\zbar\phi(w)} 
               = -~\frac{1}{\zbar\phi'(w_j)}~,
\end{equation}
where the last equality comes from the fact that $1-\zbar\phi(w_j)=0$.

Here is the result produced by (\ref{second_cphistar_formula}) and (\ref{rescalcn}): 
\begin{equation}\label{intermediate_formula}
    \cphistar f(z) ~=~ \frac{f(0)}{1-\overline{\phi(\infty)}z} 
    ~+~ \frac{1}{z} \, \sum_{j=1}^d (\overline{w_j\phi'(w_j)})^{-1} f(\wjbarinv)~,
\end{equation}
where $\seq{w_j}_{j=1}^d = \seq{1/\overline{\sigma_j(z)}}_{j=1}^d=\phi^{-1}(\{1/\zbar\})$.

Equation  \eqref{intermediate_formula} easily produces the promised formula for $\cphistar f(z)$. We require the auxiliary map 
$$
    \phibar(w) \eqdef \overline{\phi(\wbar)}~
    $$
as well as the observation that
\begin{equation}\label{der_ob}
\phibar'(w) = \overline{\phi'(\bar{w})}.
\end{equation}
 For each $w\in V$ and $j\in \{1, 2, \ldots, d\}$, we have $\phie(\sigma_j(w)) = w$, which may be rewritten
\begin{equation}\label{alt_der}
\phibar(1/\sigma_j(w)) = 1/w \ \text{for}\  w\in V\setminus\{0\}.
\end{equation}
 Differentiating both sides of \eqref{alt_der}, using \eqref{der_ob}, and substituting $z$ for $w$, we obtain for each $j$,
$$
\overline{\phi'(w_j)} = \frac{\sigma_j(z)^2}{\sigma_j'(z) z^2},
$$
which, together with $\sigma_j(z) = \wjbarinv$ and \eqref{intermediate_formula}, yields formulas \eqref{HMR_formula} and \eqref{HMR_formula_infty}.  

Note that our proof does not need to treat the case $\phinfty=\infty$ separately, since this is just the case $\alpha=0$ of (\ref{first_pfe}). It only remains to note that in this case the value $z=1/\phinftybar =0$ is now allowed, this issue having been addressed in \S\ref{discussion}  \hfill $\Box$

We can produce a more compact formula for $\cphistar$ by applying Theorem~\ref{maintheorem} and the observation that whenever $z\in \U$ is a regular value of $\phi_e$ and $\sigma_j$ is a branch of $\phie^{-1}$ defined near $z$, then $\phie'\circ\sigma_j = 1/\sigma_j'$ . Here is the the result.

\begin{corollary} \label{favorite_formula_cor}
Let $\phi\in\ratu$ have degree $d>0$. Suppose $z\in\U$ is a regular value of $\phie$ and $f\in\htwo$. Then:

{\em (a)} If $\phinfty\neq\infty$ and $z\neq 1/\phinftybar$, then
\begin{subequations} \label{HMR_formulae}
  \begin{equation}\label{HMR_formula_two}
        \cphistar f(z) = \frac{f(0)}{1-\phinftybar \, z} ~+~ 
                z \sum_{w\in\phieinv(\{z\})} \frac{f(w)}{w\phie'(w)}~.
   \end{equation}
{\em (b)} If $\phinfty=\infty$ then
    \begin{equation} \label{HMR_formula_infty_two}
        \cphistar f(z) ~=~ \left\{ 
             \begin{array}{cll}
               {\displaystyle z \sum_{w\in\phieinv(\{z\})} \frac{f(w)}{w\phie'(w)}}
                     &\quad \mbox{if} \quad & z\neq 0 \\
                     & & \\
                f(0) & \quad \mbox{if}  \quad  & z=0  
                               \end{array} \right.
    \end{equation}
    \end{subequations}
\end{corollary}

\mysubsection{Blaschke products.} \label{first_consequences}
Equations \eqref{HMR_formulae} have historical precedent in work of McDonald \cite{McD}, who considered the special case of $\phi$ a finite Blaschke product, i.e., a rational function of the form
\begin{equation} \label{blaschke_equation}
   \phi(z) = \omega\prod_{j=1}^d \frac{a_j-z}{1-\overline{a_j}\,z}~,
\end{equation}
where $\omega\in\bdu$ and $a_j\in\U$ for $1\le j\le d$.
In this case $\phie = \phi$ (since $\phi(\bdu)=\bdu$), whereupon Corollary \ref{favorite_formula_cor} yields formulae that are close in spirit to McDonald's result, but with right-hand sides that are now more explicitly computable. 

\begin{example}\label{blaschke_example} 
Following McDonald \cite{McD}, we consider Blaschke products of degree 2 with $\phi(0)=0$, i.e., with $d=2$ and $a_1=0$ in (\ref{blaschke_equation}). To simplify matters further let's take,  in that equation,  $\omega=1$ and $a_2=1/2$, so our map becomes
$$
    \phi(z) = z \, \frac{1-2z}{2-z}~.
$$
Since $\phi(\infty)=\infty$ we will compute $\cphistar$ using equation (\ref{HMR_formula_infty}), which we write formally as: 
\begin{equation}\label{blaschke_adjoint}
    \cphistar = M_{g_1} C_{\sigma_1} + M_{g_2} C_{\sigma_2}
\end{equation}
where $g_j(z) = z\sigma_j'(z)/\sigma_j(z)$, and now, because $\phie=\phi$, the $\sigma_j$'s are branches of $\phinv$.     We find
$$
    \sigma_1(z) = \frac{1+z+\sqrt{\Delta(z)}}{4} \qquad \mbox{and} \qquad
    \sigma_2(z) = \frac{1+z-\sqrt{\Delta(z)}}{4}~,
$$
where $\Delta(z) = z^2-14z+1$. If we take the square root function in these formulas to have its branch cut along the positive real axis, then $z\mapsto \sqrt{\Delta(z)}$ is holomorphic on the complex plane cut by the rays $(-\infty, 7-4\sqrt{3}]$ and $[7+4\sqrt{3},\infty)$, and is discontinuous at each point of the cut. In particular, each of the $\sigma_j$'s is holomorphic on the slit disc $\U\backslash (-1,7-4\sqrt{3}]$, and discontinuous at each point of the removed interval. So for this particular example, none of the elements on the right-hand side of (\ref{blaschke_adjoint}) has anything more than formal meaning. As we will see in the next section, this happens for {\em every\/} finite Blaschke product of degree $\ge2$. For the moment let's just notice that $7+4\sqrt{3}$ is a critical value of $\phi=\phie$ that lies in $\Ue$.
\end{example}

\section{Variants}  \label{variants}

In this section we discuss some variants of the Hammond-Moorhouse-Robbins formulae (\ref{HMR_formula})--(\ref{HMR_formula_infty}). 
In case $\phi(\infty)\in\Ue$ and $\phi(\infty)\ne \infty$, there is, as we discussed in \S\ref{discussion}, the phenomenon of pole cancellation occurring on the right-hand side of \eqref{HMR_formula}. 

\mysubsection{Finessing the ``pole problem.''} We focus, for the moment, on the subspace $\htwo_0$ consisting of functions $f\in\htwo$ with $f(0)=0$. For $f\in\htwo_0$ the first term on the right-hand side of (\ref{HMR_formula}) disappears, so with $\omegaphi$ defined, for all $z\in\U\backslash\{1/\phinftybar\}$ that are regular values of $\phie$, by
\begin{subequations}
\begin{equation}\label{omega_def}
   \omegaphi f(z) ~\eqdef  ~\sum_{j=1}^d \frac{\sigma_j'(z)}{\sigma_j(z)} f(\sigma_j(z))~= \sum_{w\in\phieinv(\{z\})} \frac{f(w)}{w\phieprime(w)}
\end{equation}
 we can write  $\cphistar f(z) = z \, \omegaphi f(z)$ (for such $z$). We know from \S\ref{much_ado} that $\cphistar f(0) = f(0) =0$, so the function $z\goesto\cphistar f(z)/z$ extends $\omegaphi f$ holomorphically to all of $\U$. {\em We will continue to use the notation $\omegaphi f$ to denote this extension.}

 Of particular interest, for $f\in\htwo_0$, is the value of our newly extended $\omegaphi f$ at $z_0 := 1/\phinftybar = \phie(0)$ when that troublesome point belongs to $\U$ and is a regular value of $\phie$.  This situation has been discussed in the third paragraph of \S\ref{discussion}, where it was observed that  $\sigma_j(z_0)=0$ for {\em exactly one\/} index $j$ which,  without loss of generality, we now take to have the value $1$. Then since $f(0)=0$ we have
   $$
       \lim_{z\goesto z_0} \frac{f(\sigma_1(z))}{\sigma_1(z)} = f'(0),
   $$ 
from which the first equation of \eqref{omega_def}, when extended to $z=z_0$, becomes:
  \begin{equation} \label{phinfty_eqn}
     \omegaphi f(z_0) = 
            \sigma_1'(z_0)f'(0) + \sum_{j=2}^d
                   \frac{\sigma_j'(z_0)}{\sigma_j(z_0)} f(\sigma_j(z_0))
  \end{equation}
  \end{subequations}
with a similar extension, which we leave to the reader for the second equation of \eqref{omega_def}. 

The next result summarizes the discussion of the last two paragraphs, using the {\em backward shift} operator introduced in \S\ref{backward_shift} to fold formulae \eqref{omega_def} and \eqref{phinfty_eqn} into a single expression.

\begin{proposition} \label{backward_prop}
 If $\phi\in\ratu$  then: 
 
 {\rm (a)}  $\omegaphi$ is a bounded linear operator from $\htwo_
0$ into $\htwo$.
 
 {\rm (b)} For each $f\in\htwo_0$ and each   $z\in\U$ that is a regular value of $\phie$,
   $$
       \omegaphi f(z) ~=~ \sum_{j=1}^d \sigma_j'(z) (Bf)(\sigma_j(z))~.
   $$
\end{proposition}
Using Proposition \ref{backward_prop} we obtain an alternative formula for $\cphistar$ in which both terms on the right-hand side represent bounded operators on $\htwo$.

\begin{corollary}\label{omegaphi_corollary}
   If $\phi\in\ratu$ then for each $f\in\htwo$:
  \begin{subequations}
  \begin{eqnarray} \label{simpler_formula_two}
     \cphistar f(z) 
          &=&~\frac{f(0)}{1-\overline{\phi(0)}\,z} 
                   + z~ \, \omegaphi[f-f(0)](z) 
                   \qquad(z\in\U)
                   \label{first_omega_formula}\\
          &~& \nonumber \\
        &=&~\frac{f(0)}{1-\overline{\phi(0)}\,z} 
                  +z\,  \sum_{j=1}^d \sigma_j'(z) (Bf)(\sigma_j(z))
                  \quad (z\in\reg(\phie)\cap\U)
                           \label{second_omega_formula}                  
  \end{eqnarray}
  \end{subequations}
  the functions $\seq{\sigma_j}_{j=1}^d$  in \eqref{second_omega_formula} being the $d$ distinct branches of $\phieinv$ defined (thanks to the $\phie$-regularity of the value $z$) in a neighborhood of $z$.
\end{corollary}
 
\begin{proof}
Upon replacing $f$ in (\ref{crucial_observation}) by the constant function $1$ we obtain:
\begin{equation} \label{cphistar_one}
    \cphistar 1(z)~ =~ \bilin{1}{K_z\circ\phi} ~=~ \overline{K_z(\phi(0))} ~=~
    \frac{1}{1-\overline{\phi(0)}\,z} \qquad (z\in\U).
\end{equation}
 Then (\ref{simpler_formula_two}) follows upon writing $f\in\htwo$ as $f=f(0)+[f-f(0)]$, using the linearity of the operator $\cphistar$,  (\ref{cphistar_one}), and the definition of $\omegaphi$.  Formula \eqref{second_omega_formula} follows from \eqref{first_omega_formula} by Proposition \ref{backward_prop}(b).
\end{proof}

\mysubsection{An amusing identity.} Upon comparing equations (\ref{HMR_formulas}), or equivalently \eqref{HMR_formulae}, with (\ref{second_omega_formula}) we find that for each regular value of $\phie$ lying in $\U\setminus\{1/\phinftybar\}$
$$
    \frac{1}{1-\overline{\phi(0)}\,z} - \frac{1}{1-\phinftybar\,z} 
            ~=~ z\,\sum_{j=1}^d\frac{\sigma_j'(z)}{\sigma_j(z)}
            ~=~ z \sum_{w\in\phieinv(\{z\})} \frac{1}{w\phie'(w)}~.
$$
This identity can be proved directly by, for example, observing that the contour integral
$$
    I(r) ~:=~ \frac{1}{2\pi i} \int_{|w|=r} \frac{dw}{w(1-z/\phie(w))}
$$
has the same value for each $r>1$, and evaluating it: (i) by residues, and (ii) by letting $r\goesto\infty$.
 
\mysubsection{Wishful thinking.} \label{wishful_thinking} 

  Proceeding in the the spirit of (\ref{blaschke_adjoint}), we can regard formula \eqref{second_omega_formula} to be a formal representation of $\cphistar$ as (modulo a rank-one operator) a sum of ``operators'' of the form $M_{z\sigma_j'} \, C_{\sigma_j } \, B$, each summand being reminiscent of Cowen's original formula for  the linear-fractional case. 
More precisely,  in equation \eqref{second_omega_formula} define the (legitimate)  rank-one operator $\Lambda_0$ on $\htwo$ by:
\begin{subequations} \label{rank_ones}
\begin{equation} \label{rank_one}
 (\Lambda_0 f)(z) 
             := \frac{f(0)}{1-\overline{\phi(0)} \, z} \qquad (f\in\htwo, z\in\U).
\end{equation} 
Then \eqref{second_omega_formula} can be rewritten, at least formally, as
\begin{equation} \label{wishful_two}
      \cphistar = \Lambda_0
                   ~+~ \sum_{j=1}^d M_{h_j} \, C_{\sigma_j} \, B,
\end{equation}
where 
   \begin{equation}
         h_j(z):= z\sigma_j'(z) \qquad (z\in\U),
   \end{equation}
   \end{subequations}
 and we make no immediate claims about the ``global legitimacy'' of the multiplication and composition operators that occur in the summation on the right-hand side of the equation.

In similar fashion, if we define 
\begin{subequations}\label{lambda_infty_eqns}
\begin{equation} \label{lambda_infty_eqn}
(\Lambda_\infty f)(z)  := \frac{f(0)}{1-\phinftybar\,z}
\end{equation}
then equation (\ref{HMR_formula}) can be rewritten formally as:
\begin{equation}\label{wishful_three}
   \cphistar = \Lambda_\infty + \sum_{j=1}^d M_{g_j}C_{\sigma_j} ,    
\end{equation}
where 
\begin{equation}
g_j(z) := \frac{z\sigma_j'(z)}{\sigma_j(z)}~,
\end{equation}
\end{subequations}
and now the legitimacy of even $\Lambda_\infty$ as an operator on $\htwo$ is an issue.

Of course our ``wishful thinking'' formulae, when applied to any $f\in\htwo$, are valid for all but finitely many $z\in\U$---that has been the point of our work up to now. In the next section we address the question of the ``global legitimacy'' of the formulae in the sense of operators on $\htwo$.

\section{Legitimate operator equations} \label{regularity} 
Let us agree to call equation (\ref{wishful_two})  a {\em legitimate operator equation\/} provided that all the multiplication and composition operators that occur on the right-hand side are bounded operators on $\htwo$, i.e., provided that the branches $\sigma_j$ of $\phieinv$ are all holomorphic on $\U$ (in which case they are automatically {\em selfmaps\/} of $\U$), and the functions $h_j$ are both holomorphic and bounded on $\U$. We make a similar definition for equation (\ref{wishful_three}), where now we must insist that $\phi(\infty)\in\U$ in order to make $\Lambda_\infty$ a bounded operator on $\htwo$. 
We saw in the Introduction (with the example $\phi(z)= z^2$) and in Example \ref{blaschke_example} above that our ``wishful thinking'' formulae (\ref{wishful_two}) and (\ref{wishful_three}) need not always be legitimate operator equations.  Here are two further examples, illustrating what can go right and what else can go wrong when we examine the legitimacy of our formulas for $\cphistar$. 

\begin{example} \label{str_outreg_example} 
Let 
$
     \phi(z) = 1/(3-z-z^2),
$
so that
$$
     \phie(w)~:=~\rho(\phi(\rho(w))) ~=~ 3 - \frac{1}{w} - \frac{1}{w^2}~.
$$
Upon solving the equation $z=\phie(w)$ for $w$ we find the two right-inverses of $\phieinv$ in the form
   $$
       \sigma_1(z) = \frac{1+\sqrt{13-4z}}{2(3-z)} 
       \qquad \mbox{and} \qquad
          \sigma_2(z) = \frac{1-\sqrt{13-4z}}{2(3-z)}~,
   $$
 where ``$\sqrt{~~}$\,'' denotes the principal branch of the square root function. Thus $\sigma_1$ and $\sigma_2$ are both defined and holomorphic on the slit plane $\C\backslash [\frac{13}{4},\infty)$ except for an isolated 
 singularity at $z=3$, which is removable for $\sigma_2$ and is a simple pole for $\sigma_1$.  What's important for our purposes is that $\sigma_1$ and $\sigma_2$ are nonzero and holomorphic on the closed unit disc.  Because they are necessarily self-maps of $\U$, they induce composition operators on $H^2$.    Since $\phi(\infty) = 0$ we see from Theorem \ref{maintheorem} that for each $f\in\htwo$ and $z\in \U\cap\reg(\phie)$:
\begin{equation}\label{sr_example}
     \cphistar f(z) = f(0) + g_1(z) f(\sigma_1(z))+g_2(z)f(\sigma_2(z)),
 \end{equation}
 where $g_j(z) = z\sigma_j'(z)/\sigma_j(z)$.   Since the $\sigma$'s are holomorphic and nonzero in a neighborhood of the closed unit disc, the $g$'s are bounded  and holomorphic in $\U$.  Because both sides of \eqref{sr_example} represent holomorphic functions on $\U$, equality holds in \eqref{sr_example} for all $z\in \U$; hence, 
 \begin{equation}\label{example_one_eqn}
     \cphistar = \Lambda_\infty + M_{g_1}C_{\sigma_1} +M_{g_2}C_{\sigma_2}
 \end{equation}
 where, because $\phi(\infty) = 0$,  $\Lambda_\infty$ is the rank one operator that takes $f\in\htwo$ to its value at the origin. So here, by contrast with the previous example, all the operators on the right-hand side of (\ref{example_one_eqn}) are defined and bounded on $\htwo$, i.e., (\ref{example_one_eqn}) is a legitimate operator equation.
  
 Observe, for future reference, that in this example the critical points of $\phi$ are $\infty$ and $-1/2$, hence the critical values are $0=\phi(\infty)$ and $4/13 = \phi(-1/2)$, both of which lie in $\U$.  \hfill $\Box$
\end{example}

The next example illustrates a situation intermediate between the behavior of Example \ref{str_outreg_example} (``good'' ) and that of the Blaschke product in Example \ref{blaschke_example} (``bad''). In this one the critical values of $\phi$ lie in $\U$, with the result that the branches $\sigma_1$ and $\sigma_2$ of $\phieinv$ are holomorphic on $\U$, but unfortunately $\phinfty\in\bdu$, which we will see causes mischief elsewhere.

\begin{example} \label{cp_bdu_example}
Let 
   $
       \phi(z)=z^2/(3-z-z^2).
   $ 
One quickly finds branches of $\phie$ in the form 
  $$
      \sigma_1(z) = \frac{1+\sqrt{\Delta(z)}}{6} \quad \mbox{and}\quad
      \sigma_2(z) = \frac{1-\sqrt{\Delta(z)}}{6}~,
 $$
 where $\Delta(z) = 13 + 12z$ and, as before, we can take ``$\sqrt{~~}$\,'' to denote the principal value of the square root. Thus both $\sigma_1$ and $\sigma_2$ are holomorphic on the slit plane $\C\backslash (-\infty,-\frac{13}{12}]$, which again contains the closed unit disc.  But now Theorem \ref{maintheorem} tells us that
 \begin{subequations} 
 \begin{equation}\label{example_two_eqn}
     \cphistar f(z) = \frac{f(0)}{1+z}
                        ~+~g_1(z)f(\sigma_1(z)) + g_2(z)f(\sigma_2(z))~,
 \end{equation}
so even though the composition operators $C_{\sigma_j}$ are both defined for $j=1,2$, at least one of the ``weights'' $g_j$ must have a pole at $-1$ to balance, when $f(0)\neq 0$, the one provided by the first term on the right. One calculates that
$$
      g_1(z) = \frac{6z}{\sqrt{\Delta(z)}(1+\sqrt{\Delta(z)})}
      \qquad\mbox{and}\qquad
     g_2(z) =  \frac{-6z}{\sqrt{\Delta(z)}(1-\sqrt{\Delta(z)})}~,
$$
so that that $g_1$ is bounded on $\U$, while $g_2$ has required the pole at $-1$. 
Thus in  the ``wishful thinking'' formula
\begin{equation}\label{partially_legit_formula}
    \cphistar ~=~ \Lambda_\infty f ~+~ 
                     M_{g_1}C_{\sigma_1} ~+~ M_{g_2}C_{\sigma_2}
\end{equation}
\end{subequations}
derived from (\ref{example_two_eqn}), all the operators make sense when acting on the space $H(\U)$ of all functions holomorphic on $\U$, but on $\htwo$, even though $M_{g_1}$ and both composition operators  act boundedly,  $M_{g_2}$ and $\Lambda_\infty$ do not. In short, we may regard (\ref{partially_legit_formula}) as a ``legitimate operator equation'' on $H(\U)$, but not on $\htwo$. 

On the other hand,  note that $h_j(z) = z\sigma_j'(z)$ is bounded and analytic on $\U$ for $j=1$ and $j=2$.  Thus, we can regard (\ref{wishful_two}) as a legitimate operator equation for $\phi$, $\sigma_1$, and $\sigma_2$ as above.   \hfill $\Box$
\end{example}

The examples given so far show that the location both of critical values of $\phi$ and of $\phinfty$ are crucial to the way in which formulas (\ref{wishful_two}) and (\ref{wishful_three}) can be interpreted. Let us  now make this official. 

\begin{definitions} 
{\em (a)} We call a map $\phi\in\ratu$ {\em outer regular\/} whenever its critical values all lie in $\U$, i.e., whenever each point of $\overline{\Ue}:=\seq{|z|\ge 1}$ is a regular value of $\phi$. 

{\em (b)} If, in addition, $\phinfty\in\U$ we say that $\phi$ is {\em strongly\/} outer regular. 
\end{definitions}

Thus the mapping of Example \ref{str_outreg_example} is strongly outer regular, while that of Example \ref{cp_bdu_example}  is outer regular, but, because $\phinfty$ lies on $\bdu$,  not strongly so. The Blaschke product of Example \ref{blaschke_example} has a critical value in $\Ue$, so it is not outer regular. This, in fact, is a feature of {\em every\/} finite Blaschke product of degree $\ge2$.

\begin{proposition} \label{blaschke_nonreg} 
Finite Blaschke products of degree $\ge 2$ are {\em never\/} outer regular.
\end{proposition}

\begin{proof}
  A finite Blaschke product $\phi$ of order $>1$ must have a critical point somewhere in the finite plane (the derivative must vanish somewhere). This cannot occur on $\bdu$, else $\phi$ could not take the unit disc into itself.  Because the unit circle is taken onto itself, the reflection in the unit circle of any critical point is again a critical point, so there is a critical point in $\Ue$. Since Blaschke products are also selfmaps of $\Ue$ there is a critical value in $\Ue$.
  \end{proof}

Here is our main result on outer regularity.

\begin{theorem} \label{global_thm}
Suppose $\phi\in\ratu$ has degree $d$.

{\em (a)} If $\phi$ is outer regular then  {\em (\ref{wishful_two})} is a legitimate operator equation.
   
{\em (b)}  If, in addition, $\phi$ is {\em strongly\/} outer regular, then   {\em (\ref{wishful_three})} is a legitimate operator equation.
   \end{theorem}

\begin{proof}
(a) The hypothesis of outer regularity on $\phi$ means that  $r\Ue\subset \reg(\phi)$ for some $0<r<1$, hence $(1/r)\U \subset \reg(\phie)$. Thus Corollary \ref{sc_continuation} guarantees that $\phieinv$ has $d$ distinct branches $\sigmaseq$ on $(1/r)\U$, each of which, as we noted in \S\ref{ratselfmap_intro}, maps $\U$ into itself. Since each of these branches is holomorphic in a neighborhood of $\ubar$, each has a derivative bounded on $\U$, hence the functions $h_j$ of formula (\ref{wishful_two}) are bounded on $\U$, and induce bounded multiplication operators on $\htwo$.

(b) The boundedness of $\Lambda_\infty$ on $\htwo$ is a consequence of the fact $\phi(\infty)\in\U$ ({\em strong\/} outer regularity). The boundedness of the functions $g_j$ will follow from that of the derivatives of the $\sigma_j$'s once we establish that each $\sigma_j$ is bounded away from zero on $\U$. Now strong outer regularity guarantees that $\phie(0)\in\Ue$. Fix $j$ and suppose, for the sake of contradiction that $\sigma_j(z)=0$ for some $z\in\ubar$. Then, recalling that $\sigma_j$ is a branch of $\phieinv$ on a neighborhood of the {\em closure\/} of the unit disc,  $z=\phie(\sigma_j(z)) = \phie(0)\in\Ue$, a contradiction. Thus $\sigma_j$ does not vanish at any point of $\ubar$, and so it is bounded away from zero on $\U$.
\end{proof}

\section{Adjoints modulo compacts}  \label{whendoes}
In this section we set out conditions under which adjoint composition operators are compact perturbations of simpler ones, e.g., weighted composition operators.  The following example, which occurred in \S\ref{regularity} as the poster child for the concept of strong outer regularity, motivates our work.

\begin{example} \label{pauls_example}
  {\em  If
   $
       \phi(z)=1/(3-z-z^2) 
   $   
    then $\cphistar$ is a compact perturbation of a weighted composition operator.}
 \end{example}  
 
\begin{proof}  Recall from Example \ref{str_outreg_example} the branches $\sigma_1$ and $\sigma_2$ of $\phieinv$, both of which are holomorphic on a neighborhood of the closed unit disc.  Note that $\sigma_1(1)=1$, but that $\sigma_2$ takes the closed unit disc into $\U$. Thus  $C_{\sigma_2}$ is compact on $\htwo$ \cite[\S2.2, page 23]{JS}  while $C_{\sigma_1}$ is not \cite[\S 4.1, ``Proposition'' on page 56]{JS}, so it follows from (\ref{example_one_eqn}) that $\cphistar$ is a compact perturbation of the weighted composition operator $M_{g_1}C_{\sigma_1}$. 
\end{proof}

Here is the general result motivated by this example.

\begin{corollary} \label{calkin_cor}
Suppose $\phi\in\ratu$ is outer regular and maps exactly one point of $\bdu$ to $\bdu$. Then:

{\em (a)} $\cphistar$ is a compact perturbation of an operator of the form $M_hC_\sigma B$ where $\sigma$ is a holomorphic selfmap of $\U$ that is a branch of $\phieinv$, the function $h:z\goesto z\goesto z\sigma'(z)$ is holomorphic and bounded on $\U$, and $B$ is the backward shift on $\htwo$.  

{\em (b)} If, in addition, $\phi$ is {\em strongly\/} outer regular, then $\cphistar$ is a compact perturbation of the weighted composition operator $M_g C_\sigma$ with $\sigma$ as in part (a) and $g:z \goesto z\sigma'(z)/\sigma(z)$ a bounded holomorphic function on $\U$.
\end{corollary}

\begin{proof}
(a) The hypothesis on the boundary behavior of $\phi$ insures that $\phieinv(\bdu)$ contains exactly one point of $\bdu$.   By \S\ref{ratselfmap_intro} and Proposition \ref{disjointness_prop} the branches $\sigmaseq$ map the closed unit disc $\ubar$ to pairwise disjoint closed subsets of $\ubar$, so the boundary behavior imposed on $\phie$ guarantees that exactly one of the maps $\sigma_j$, say $\sigma_1$, takes a point of $\bdu$ to $\bdu$, with the other $\sigma_j$'s  taking $\bdu$ into $\U$. Thus, as in the discussion of Example \ref{pauls_example}, the composition operators $C_{\sigma_j}$ are compact for $2\le j\le d$, while $C_{\sigma_1}$ is {\em not\/} compact. It follows that each of the operators $M_{z\sigma_j'}C_{\sigma_j}B$ is compact ($2\le j \le d$) so, upon noting that the rank-one operator $\Lambda_0$ is also compact, part (a) of Theorem \ref{global_thm} provides the desired representation for $\cphistar$, with $\sigma:=\sigma_1$. 

(b) Use part (b) of Theorem \ref{global_thm} exactly as above, noting that, again as in Example \ref{pauls_example} the strong outer regularity assumed for $\phi$ guarantees the boundedness, hence compactness, of the rank-one operator $\Lambda_\infty$.
\end{proof}

As an application of Corollary \ref{calkin_cor}, consider the following generalization of Example \ref{pauls_example}.

\begin{example} {\em
For $d=2, 3, \, \ldots \,$ let 
  $$
     \phi(z) = \frac{1}{(d+1)-z-z^2 \, \cdots \, - z^d} ~.
  $$ 
Then $\cphistar$ is a  compact perturbation of a weighted composition operator.  
}
\end{example}

\begin{proof} One checks easily that $\phi$ maps the closed unit disc into $\U\cup\{1\}$,  taking 1 to 1 and the rest of the unit circle into $\U$. Thus,
in view of part (b) of Corollary \ref{calkin_cor}, we need only show that $\phi$ is strongly outer regular. Since $\phinfty=0\in\U$, outer regularity is the only issue, i.e. we must show that all the critical values of $\phi$ lie in $\U$. Note first that $\infty$ is a critical point, with the corresponding critical value $\phinfty=0\in\U$.   As for the critical points in the finite plane: these are the critical points of the denominator---the points $z$ at which  the derivative of the denominator vanishes; equivalently, the points $z$ for which:

\begin{equation}
    f_d(z) := 1+2z +3z^2 + \, \cdots \, +dz^{d-1} = 0~.
\end{equation}
It will be easier to deal with
$$
   g_d(z) := (1-z)f_d(z) = 1 + z + z^2 + \,\cdots\, z^{d-1} - dz^d,
$$
which adds to the zeros of $f_d$ an additional one at $z=1$. A simple exercise involving Rouch\'e's theorem shows that $g_d$ has, for each $R>1$, $d$ zeros in the disc $R\U$, hence $d$ zeros in $\ubar$, and therefore $f_d$ has $d-1$ zeros  (counting multiplicity) in $\ubar$. To complete the proof note that $\phi$ takes the closed disc into the open disc with the exception of the point $1$. Since $f_d(1)\neq 0$ we see that all critical values of $\phi$ belong to $\U$.
\end{proof}

It's not difficult to construct other examples of outer regular and strongly outer regular mappings; e.g., $\phi(z) = z/(a-z^n)$ is a self-map of $\U$ whenever $|a| \ge 2$ and is strongly outer regular whenever the integer $n$ exceeds $1$. 


\end{document}